\documentclass{svmultmod}

\usepackage{type1cm}
\usepackage{graphicx}
\usepackage{newtxtext}
\usepackage[varvw]{newtxmath}

\newcommand{\NN}{\mathbb{N}}
\newcommand{\ZZ}{\mathbb{Z}}
\newcommand{\RR}{\mathbb{R}}
\newcommand{\CC}{\mathbb{C}}
\newcommand{\HH}{\mathbb{H}}
\newcommand{\gammafunc}{\mathrm{\Gamma}}
\newcommand{\mc}{\mathcal}
\DeclareMathOperator{\PSL}{PSL}
\DeclareMathOperator{\SL}{SL}
\DeclareMathOperator{\GL}{GL}
\DeclareMathOperator{\Unit}{U}
\DeclareMathOperator{\Rea}{Re}
\DeclareMathOperator{\Hei}{Hei}

\theoremstyle{plain}
\newtheorem{thm}{Theorem}[section]

\usepackage[style=rwtmath, dashed=true]{biblatex}
\begin{document}

\title*{Some aspects of the spectral theory with twisting representations}
\titlerunning{Anke Pohl: Aspects of spectral theory with twists}
\author{Anke Pohl}
\institute{Anke Pohl \at University of Bremen, Department~3: Mathematics and Computer Science, Institute for Dynamical Systems, Bibliothekstr.~5,  28359~Bremen, Germany, \email{apohl@uni-bremen.de}}

\maketitle

\abstract{In recent years, significant advancements have been made in understanding the spectral theory of hyperbolic spaces, particularly in the context of twisting representations, both unitary as well as non-unitary. We survey some results that we obtained within the framework of the SPP~2026 ``Geometry at Infinity'' as well as some closely related results.}

\keywords{spectral theory, hyperbolic spaces, Laplace eigenfunctions, Jacobi Maass forms, unitary twists, non-unitary twists, non-expanding cusp monodromy, Selberg zeta function, divisor, transfer operator, spectral interpretation of zeros}

\medskip

\noindent
\textbf{2020 Mathematics Subject Classification.} {Primary: 11M36, 58J50; Secondary: 30F35, 11F12, 30F40}

\section{Introduction}
\label{sec:intro}

The spectral theory of Riemannian manifolds or, more generally, Riemannian orbifolds is fundamental for many research areas within mathematics, physics as well as for various applications. Accordingly, spectral theory features many different research directions. The main focus in this article is on the study of the relation between geometric or dynamical entities on the one hand and spectral entities on the other hand, in particular between geodesics and their length spectrum on the dynamical side and Laplace eigenfunctions and resonances on the spectral side. Understanding geodesics as the classical mechanical description of (friction-less, undisturbed) motions on the given Riemannian orbifold, and Laplace eigenfunctions as the relevant building blocks of the corresponding quantum mechanical description at the subatomic level, results on the interdependence of properties of geodesics and Laplace eigenfunctions can often be seen as mathematical instances of Bohr's correspondence principle between classical mechanics and quantum mechanics in physics.

In this article we will mainly consider hyperbolic orbisurfaces, thus spaces of the form $X = \Gamma\backslash\HH$ with $\HH$ being the hyperbolic plane and $\Gamma$ the fundamental group of~$X$. We may and shall use the upper half plane model for~$\HH$ and identify $\Gamma$ with a Fuchsian group, i.e., a discrete subgroup of~$\PSL_2(\RR)$, acting via fractional linear transformations on~$\HH$. We emphasize that $\Gamma$ may contain elliptic elements, in which case $X$ is a genuine orbifold, not a manifold.

Classically, $\Delta_X$, the Laplacian of~$X$, is considered as an operator on some suitable space of functions $X\to\CC$. However, any function $f\colon X\to\CC$ can be identified with a $\Gamma$-invariant function $u\colon \HH\to\CC$ via
\[
 u(z) = f([z])
\]
for any $z\in\HH$ with $[z]=\Gamma.z$ denoting the $\Gamma$-orbit of~$z$ as an element of~$X=\Gamma\backslash\HH$. Considering $\Delta_X$ as acting on suitable spaces of functions $u\colon\HH\to\CC$ obeying the $\Gamma$-invariance
\[
 u(g.z) = u(z)\qquad\text{for all $g\in\Gamma$, $z\in\HH$}\,,
\]
immediately opens up the possibility of generalizing this setup to an equivariance situation instead of a pure invariance situation. To be more precise, let $V$ be a finite-dimensional complex vector space and let $\chi\colon\Gamma\to\GL(V)$ be a representation. The Laplacian~$\Delta_X$ obviously acts on vector-valued functions $u\colon \HH\to V$ of suitable regularity. We will consider $\Delta_X$ as acting on suitable spaces of such functions obeying the equivariance property
\begin{equation}\label{eq:equivariance}
 u(g.z) = \chi(g)u(z)\qquad\text{for all $g\in\Gamma$, $z\in\HH$}\,.
\end{equation}
We typically refer to the representation~$\chi$ in this setup as a \emph{twist}.

We may equivalently model this setup using the flat vector orbibundle with typical fiber~$V$ associated to~$\chi$. This is determined by the diagonal action of~$\Gamma$ on~$\HH\times V$, i.e.,
\[
g.(z,v) = (g.z, \chi(g)v)
\]
for $g\in\Gamma$, $(z,v)\in\HH\times V$. The vector orbibundle is the space of $\Gamma$-orbits on~$\HH\times V$,
\[
\HH\times_\chi V \coloneqq \Gamma\backslash (\HH\times V)\,,
\]
with the canonical projection on the $\HH$-component as bundle projection. The \mbox{(Bochner--)} Laplacian acts on sections of~$\HH\times_\chi V$ of suitable regularity. Identifying these sections with functions satisfying~\eqref{eq:equivariance} gives the equivalence of the two models.

Already Selberg~\cite{Selberg, Selberg_est_fourier} promoted the idea to investigate twisted situations alongside untwisted one and to study their interrelations in order to understand the properties of their entities. Even though research initially was restricted to unitary twists, the same motivations apply to non-unitary twists as well. We list a few reasons for studying twisted situations.

\begin{itemize}
\item Considering the classical, untwisted entities as the primary object of interest, e.g., Laplace eigenfunctions~$u$ with $\Delta_X u = \lambda u$ and $u(g.z)=u(z)$, we may consider \eqref{eq:equivariance} as a perturbation of these objects. More precisely, instead of studying a single twisted situation, we should rather ask for a family of twists~$\chi_t$ depending in some controlled way of a parameter~$t$ (e.g., continuously or smoothly) with $\chi_0=1$. By investigating which entities change in which way while slowly passing away from the untwisted situation $\chi_0$ along the given family, we may identify important features, invariants, rigidities or non-rigidities.
\item The entities from twisted situations are of interest on their own. To name but a few examples, they play an important role in the investigation of Casimir interaction between two conducting obstacles~\cite{Li_Strohmaier}, in studying tunneling problems~\cite{Ludewig_Rosenberger}, in questions on black holes and AdS$_3$ space-times~\cite{bytsenko2006spectral, DGG, MatsushimaMurakami}. Analogously, vector-valued modular forms (even though not Laplace eigenfunctions but nevertheless underlying the same line of thought) are of great importance in the study of the generalized moonshine conjecture~\cite{Dong_Li_Mason, Duncan_ODesky,gannon2006monstrous} and in two-dimensional conformal field theory~\cite{francesco2012conformal, Harvey_Wu}.
\item Considering several situations with or without twist simultaneously can help us to take advantage of symmetries to break down considerations into smaller pieces. Prototypical examples of this type are decomposition formulas such as the decomposition of representations into irreducible factors, which typically get inherited to Laplace eigenfunctions, or formulas of Venkov--Zograf type that allow us to express the Selberg zeta function twisted by some representation as the product of Selberg zeta functions twisted by factors of this representation. See, e.g., \cite{Venkov_Zograf, FP_szf}. In a similar vein, such considerations can allow us to  explore inherent relations between different situations to study certain questions from different points of view. Prototypical examples are again formulas of Venkov--Zograf-type that allow us to relate the Selberg zeta function of a group and its subgroups by well-chosen representations.
\end{itemize}

Situations with \emph{unitary} twisting representations have been investigated with great success for several decades by now, providing a plethora of results to which current research still adds and to which also some of our results obtained within the framework of the SPP~2026 contribute. In Section~\ref{sec:unitary} we will survey some of our recent achievements.

The investigation of situations with \emph{non-unitary} twisting representations got increased interest much more recent and many, also many rather elementary questions are still open. A (necessarily incomplete) list of results for non-unitary twists at the start of the second phase of the SPP~2026 includes~\cite{Eholzer_Skoruppa, Daughton2016,
Raji_cohom, Muehlenbruch_Raji, Deitmar_Monheim_traceformula, Deitmar_Monheim_eisenstein, Deitmar_locally_compact, FP_szf, FP_eisenstein, Saber_Sebbar, Knopp_Mason_definition, Knopp20032, Knopp_Mason_illinois, Knopp2011, Knopp_Mason2012, Kohnen_Mason, Kohnen_Martin, Kohnen, Kohnen_Mason2, Bantay_Gannon, Creutzig_Gannon, Gannon_theory, DGMS} for hyperbolic surfaces, and \cite{Mueller_STF, Spilioti2015, Fedosova_nonunitary,
Deitmar_Monheim_traceformula, Deitmar_locally_compact, Spilioti2020} for higher-dimensional spaces. In Section~\ref{sec:nonunitary} we will discuss some of our contributions.

\begin{acknowledgement}
AP's research is funded by the Deutsche Forschungsgemeinschaft (DFG, German Research Foundation) -- project no.~441868048 (Priority Program~2026 ``Geometry at Infinity''). She likes to thank Claudio Meneses and Christopher Wulff for helpful comments.
\end{acknowledgement}

\section{A first example}
\label{sec:example}

We illustrate with a small example that the passage from twists by unitary representations to twists by arbitrary representations indeed allows for new observations and results. To that end, we let $\Gamma$ be a Fuchsian group generated by a single hyperbolic element, say
\[
\Gamma = \langle a_\ell \rangle = \{ a_\ell^n : n\in\ZZ\}
\]
with
\[
a_\ell = \begin{pmatrix} e^{\ell/2} & 0 \\ 0 & e^{-\ell/2}  \end{pmatrix}\qquad\text{for some $\ell>0$}.
\]
The hyperbolic surface $X = \Gamma\backslash\HH$, commonly known as a \emph{hyperbolic cylinder}, has two funnels, no cusps and no orbifold points. It has two primitive periodic geodesics, which both trace out the same path but are time-inverted to each other. Both are of length~$\ell$. The \emph{Selberg zeta function} of~$X$ is given by the infinite product
\[
Z_X(s) \coloneqq \prod_{k=0}^\infty \left(1 - e^{-(s+k)\ell} \right)^2 \qquad\text{for $s\in\CC$\,, $\Rea s >0$\,,}
\]
and its meromorphic (indeed holomorphic) extension to all of~$\CC$. We denote the extension by~$Z_X$ as well. The zeros of~$Z_X$ coincide with the Laplace resonances of~$X$. We refer to Sections~\ref{sec:unitary} and~\ref{sec:nonunitary} for more details and take here for granted that the zeros of~$Z_X$ are of interest for the spectral theory of~$X$. A straightforward calculation shows that the zero set of~$Z_X$ is
\[
-\NN_0 + i\frac{2\pi}{\ell}\ZZ\,,
\]
each zero having order~$2$.

Now let $\chi\colon\Gamma \to \CC^\times$ be a one-dimensional representation and fix $\varphi\in\CC$ such that
\[
\chi(a_\ell) = e^{i\varphi}\,.
\]
The \emph{twisted Selberg zeta function} of~$X$ and~$\chi$ is given by the infinite product
\[
Z_{X,\chi}(s) \coloneqq \prod_{k=0}^\infty \left( 1 - e^{i\varphi - (s+k)\ell} \right)\left(1-e^{-i\varphi - (s+k)\ell}\right)\quad\text{for $\Rea s \gg 1$\,,}
\]
and its meromorphic (indeed holomorphic) continuation to all of~$\CC$. We again refer to Sections~\ref{sec:unitary} and~\ref{sec:nonunitary} for more details. The zero set of~$Z_{X,\chi}$ is
\[
-\NN_0 \pm i\frac{\varphi}{\ell} + i\frac{2\pi}{\ell}\ZZ\,,
\]
each zero is of order~$1$, unless it is represented in several ways in this set (in which case it is of order~$2$).

We now consider a one-parameter family of representations~$(\chi_t)_{t\in [0,\infty)}$ starting with the trivial representation~$\chi_0=1$, and we investigate how the zeros of the associated (twisted) Selberg zeta functions move along with the parameter~$t$. In other words, we want to understand how the zeros of the untwisted Selberg zeta function~$Z_X$ behave under perturbations provided by~$(\chi_t)_t$. By direct calculation we observe that along a \emph{unitary} family, e.g., $\chi_t(a_\ell) = e^{it}$, the zeros only change their imaginary part. Thus, using only unitary representations, a relation between different zeros of~$Z_X$ can be observed only for zeros with same real part.

However, if we allow \emph{arbitrary} representations as well, e.g., $\chi_t(a_\ell) = e^{i(2+3i)t}$, then the zeros move along much more general paths, in particular, a change in their real parts is also possible. In this way, we can observe relations between any two different zeros of~$Z_X$.

\section{Spectral theory for unitary twists, and Selberg zeta functions for hyperbolic orbisurfaces of infinite area}
\label{sec:unitary}

In this section we consider a geometrically finite hyperbolic orbisurface~$X=\Gamma\backslash\HH$ and a finite-dimensional \emph{unitary} twist~$\chi\colon \Gamma \to \Unit(V)$ of its fundamental group~$\Gamma$. We emphasize that we allow $\Gamma$ to have elliptic elements. If $X$ is of finite area, or if $\chi$ is trivial and $\Gamma$ is torsion-free, then the Selberg zeta function (defined in Section~\ref{sec:finarea}) is known to mediate between the geodesic length spectrum of~$X$ and the resonances of the (twisted) Laplacian on~$X$. With our work we contribute to extending these results to the general case of arbitrary geometrically finite Fuchsian groups~$\Gamma$ and arbitrary finite-dimensional unitary twists. To be more precise, we start, in Sections~\ref{sec:finarea} and~\ref{sec:notwist}, with a more detailed survey of the preliminary results that motivated our work. We then survey our own results in Section~\ref{sec:ours}.

\subsection{Finite area hyperbolic orbisurfaces}
\label{sec:finarea}

We suppose for the moment that $X = \Gamma\backslash\HH$ is of finite area and that $\chi$ is the trivial one-dimensional representation, for which reason we omit it from the notation. The \emph{Selberg zeta function}~$Z_X$ for~$X$ is determined by the infinite product
\begin{equation}\label{def:szf}
Z_X(s) \coloneqq \prod_{\ell\in L(X)} \prod_{k=0}^\infty \left( 1 - e^{-(s+k)\ell} \right)\,,
\end{equation}
where $L(X)$ denotes the multiset of lengths of primitive periodic geodesics on~$X$. This infinite product converges for $\Rea s > 1$, and it admits a meromorphic continuation to all of~$\CC$, which we continue to denote by~$Z_X$. See, e.g., \cite{Selberg, Venkov_book, Hejhal1, Hejhal2}.

In the domain $\{\Rea s \leq 1\}$, i.e., in the domain where the Selberg zeta function~$Z_X$ is defined by meromorphic continuation, $Z_X$ may have zeros. To understand the importance of these zeros, we consider the resolvent
\[
R_X(s) \coloneqq \left( \Delta_X - s(1-s) \right)^{-1}
\]
of the (positive) Laplacian~$\Delta_X$ on~$X$. The resolvent~$R_X$ defines an operator on~$L^2(X)$ if $\Rea s > 1/2$ and simultaneously $s(1-s)$ is not an $L^2$-Laplace eigenvalue. It admits a meromorphic continuation to all of~$\CC$ as a family of operators
\[
R_X(s) \colon L^2_{\text{comp}}(X) \to L^2_{\text{loc}}(X)
\]
with poles of finite rank (see, e.g., \cite{Mazzeo_Melrose}). These poles are the \emph{resonances} of~$X$. We emphasize that each spectral parameter of~$\Delta_X$ (considered as an operator on $L^2$-functions) is indeed a resonance, and hence every $L^2$-Laplace eigenvalue gives rise to two resonances (or twice the resonance $1/2$). The zeros of~$Z_X$ coincide (up to certain well-understood exceptions of minor magnitude) with the resonances, including their orders. See, e.g., \cite{Selberg, Venkov_book, Guillope, Hejhal1, Hejhal2}.

If $\chi\colon\Gamma\to \Unit(V)$ is now any finite-dimensional representation of~$\Gamma$, then these results generalize to the unitarily twisted situation. More precisely, the \emph{twisted Selberg zeta function}~$Z_{X,\chi}$ is then determined by the infinite product
\begin{equation}\label{def:szftwist}
Z_{X,\chi}(s) \coloneqq \prod_{[g]\in [\Gamma]_{\text{hp}}} \prod_{k=0}^\infty \det\left( 1 - \chi(g)e^{-(s+k)\ell(g)} \right)\,,
\end{equation}
where we use the canonical bijection between $\Gamma$-conjugacy classes~$[\Gamma]_{\text{hp}}$ of primitive hyperbolic elements in~$\Gamma$ and the primitive periodic geodesics on~$X$, and let $\ell(g)$ denote the (primitive) length of the periodic geodesic associated to~$g$. This infinite product is still convergent for~$\Rea s > 1$, and it continues to admit a meromorphic continuation to all of~$\CC$. See, e.g., \cite{Selberg, Venkov_book, Hejhal1, Hejhal2}. Also on the spectral side, in the same way as for the untwisted situation, the twisted Laplacian~$\Delta_{X,\chi}$ has a resolvent
\[
R_{X,\chi}(s) \colon L^2(X,E_\chi) \to L^2(X,E_\chi)\,,
\]
with $E_\chi \coloneqq \HH\times_\chi V$, being well-defined for $\Rea s > 1/2$ and $s(1-s)$ not being a Laplace eigenvalue which, considered as an operator
\[
R_{X,\chi}(s) \colon L^2_{\text{comp}}(X,E_\chi) \to L^2_{\text{loc}}(X,E_\chi)\,,
\]
extends to a family such operators that is meromorphic on all of~$\CC$. The poles of this family, i.e., the resonances of~$\Delta_{X,\chi}$, are in bijection with the zeros of~$Z_{X,\chi}$, again up to minor, well-understood exceptions. See, e.g., \cite{Selberg, Venkov_book, Hejhal1, Hejhal2, Phillips_perturb_twist, Phillips_scattering_twist}.

\subsection{Infinite area hyperbolic surfaces, no orbifold points, no twists}
\label{sec:notwist}

For hyperbolic orbisurfaces, the methodology used for \emph{finite} area spaces sometimes does not carry over to \emph{infinite} area spaces. Indeed, not all techniques used to obtain the results surveyed in Section~\ref{sec:finarea} can be applied to infinite area hyperbolic orbisurfaces. Nevertheless, some of the results could be generalized to infinite area situations, and these results were the main motivation for our work surveyed in Section~\ref{sec:ours} below.

Let now $X=\Gamma\backslash\HH$ be a hyperbolic \emph{surface} of infinite area, i.e., $\Gamma$ has no elliptic elements and hence $X$ has no orbifold singularities. The definition of a Selberg zeta function~$Z_X$ for~$X$ as in~\eqref{def:szf} still applies, which means that the infinite product converges for $\Rea s \gg 1$. More precisely, in this situation the axis of convergence is the Hausdorff dimension~$\delta$ of the limit set of~$\Gamma$. (Indeed, it also is so for cofinite Fuchsian groups as the Hausdorff dimension is then~$1$.) In~\cite{Patterson_Perry} the Selberg zeta function for hyperbolic surfaces without cusps was studied, and in~\cite{BJP} the investigations were extended to hyperbolic surfaces with cusps. The result by Borthwick, Judge and Perry in the latter work provides a meromorphic continuation of~$Z_X$ and, more importantly, a full understanding of the divisor of~$Z_X$.

\begin{thm}[{\cite[Theorem~4.1]{BJP}}]\label{thm:BJP}
Let $X$ be a non-elementary geometrically finite hyperbolic surface of infinite area. Then the Selberg zeta function admits the factorization
\[
Z_X(s) = e^{p(s)} G_\infty(s)^{-\chi_e^{\text{top}}(X)}\gammafunc\bigl(s-\tfrac12\bigr)^{n_p} \mc P_X(s)
\]
on all of~$\CC$.
\end{thm}

In Theorem~\ref{thm:BJP}, $\mc P_X$ is the Hadamard product of order~$2$ of the resonances of~$X$, i.e.,
\[
\mc P_X(s) = s^{m_X(0)} \prod_{r\in \mc R_X\setminus\{0\}} E_2\left( \frac{s}{r}\right)\,,
\]
where $\mc R_X$ denotes the multiset of resonances of~$X$, the integer~$m_X(0)$ is the order of~$0$ as a resonance, and $E_2$ is the second elementary factor. Further,
$p$ is a polynomial of degree at most~$2$, the number~$\chi_e^{\text{top}}(X)$ is the topological Euler characteristic of~$X$, the function $G_\infty$ is essentially a certain product of the Barnes $G$-function~$G$, namely
\[
G_\infty(s) = \frac{1}{(2\pi)^s} G(s) G(s+1)\,.
\]
Finally, $\gammafunc$ denotes the gamma function, and $n_p$ is the number of cusps of~$X$. As $G_\infty$ is entire with zeros only at the negative integers, and $\gammafunc$ is meromorphic without zeros, Theorem~\ref{thm:BJP} shows that $Z_X$ has a zero if and only if this zero is a resonance of~$X$, unless there is interference between resonances and poles/zeros of the factor functions. The latter situations, however, are well-understood. We refer to the textbook~\cite{Borthwick_book} for detailed discussions of these results and for further references. For finite area hyperbolic orbisurfaces, an analogous divisor formula can be deduced from~\cite{Venkov_book}.

\subsection{Infinite area with orbifold points and unitary twists}
\label{sec:ours}

Motivated by Theorem~\ref{thm:BJP}, the author of this note together with M.~Doll and partly together with K.~Fedosova obtained a generalization of Theorem~\ref{thm:BJP} to arbitrary geometrically finite hyperbolic orbisurfaces~$X=\Gamma\backslash\HH$ (i.e., allowing also elliptic elements in~$\Gamma$) and twists by arbitrary finite-dimensional unitary representations~$\chi\colon\Gamma\to\Unit(V)$. See \cite{DFP, DFP2, Doll_Pohl}.

In this situation, the infinite product in~\eqref{def:szftwist} converges for $\Rea \gg 1$, indeed for $\Rea s > \delta$ (see, e.g., \cite{FP_szf}) and hence can be used as a starting point for the definition of a twisted Selberg zeta function~$Z_{X,\chi}$. To obtain a generalization of Theorem~\ref{thm:BJP}, one needs to
\begin{itemize}
\item provide a suitable definition of a Laplacian in the presence of the combination of infinite area, orbifold points and twisting representations,
\item understand the relevant parts of its spectral theory,
\item provide a meromorphic continuation of the infinite product used for the definition of the Selberg zeta function, and
\item understand the structure of the Selberg zeta function in relation to the spectral theory of the Laplacian, in particular in the presence of the combination of elliptic elements and twisting representations.
\end{itemize}
Selberg's Lemma guarantees the existence of a finite-index subgroup, $\Lambda$, of~$\Gamma$ without elliptic elements. Furthermore, \cite{FP_szf} provides formulas of Venkov--Zograf-type for infinite-area situations as well. Nevertheless, a descent from $\Gamma$ to~$\Lambda$ combined with a twisting by the induced representation from~$\Lambda$ to~$\Gamma$ does not yield a full generalization of Theorem~\ref{thm:BJP}, not even for untwisted situations. For this reason, the contribution of elliptic elements has to be handled directly in each step of the proofs.

The Laplace operator~$\Delta_{X,\chi}$ can be defined, as already indicated in Section~\ref{sec:intro}, rather straightforwardly as an operator acting on vector-valued functions obeying the equivariance relation~\eqref{eq:equivariance}, or, equivalently, as the (Bochner-)Laplacian acting on sections of~$E_\chi = \HH\times_\chi V$.

\begin{thm}[\cite{DFP}]\label{thm:DFP}
The Laplace operator~$\Delta_{X,\chi}$ is self-adjoint on~$L^2(X, E_\chi)$. Its resolvent~$R_{X,\chi}$, considered as an operator family
\[
R_{X,\chi}(s) \colon L^2_{\text{comp}}(X,E_\chi) \to L^2_{\text{loc}}(X,E_\chi)\,,
\]
admits a meromorphic continuation to a family of such operators on all of~$\CC$.
\end{thm}

The proof technique of Theorem~\ref{thm:DFP} extends the one by~\cite{GuZw95} to twisting representations and the presence of elliptic elements. The structure in the proof that funnel ends, cusp ends, and the compact core are partly handled separately remains preserved and eventually also reflects in the divisor formula for the twisted Selberg zeta function. See Theorem~\ref{thm:DP} below.

A crucial step in obtaining a divisor formula for the twisted Selberg zeta function, is to establish a factorization formula for the twisted scattering determinant~$\tau_{X,\chi}$.

\begin{thm}[\cite{DFP2}]\label{thm:DFP2}
The twisted scattering determinant~$\tau_{X,\chi}$ admits the factorization
\[
\tau_{X,\chi}(s) = e^{p(s)}\,\frac{\mc P_{X,\chi}(1-s)}{\mc P_{X,\chi}(s)}\,\frac{\mc P_{X_f,\chi}(s)}{\mc P_{X_f,\chi}(1-s)}\,,
\]
where $\mc P_{X,\chi}$ is the Hadamard product of the resonances of~$\Delta_{X,\chi}$, the space~$X_f$ is the disjoint union of the funnel ends, $\mc P_{X_f,\chi}$ is the Hadamard product associated to the resonances of~$X_f$, and $p$ is a polynomial of degree at most~$4$.
\end{thm}

These results come together in the following result on the twisted Selberg zeta function.

\begin{thm}[\cite{Doll_Pohl}]\label{thm:DP}
For any geometrically finite hyperbolic orbisurface~$X=\Gamma\backslash\HH$ of infinite area and for any finite-dimensional unitary representation~$\chi\colon\Gamma\to\Unit(V)$, the infinite product~\eqref{def:szftwist} has a meromorphic continuation to all of~$\CC$, which we denote by~$Z_{X,\chi}$. The twisted Selberg zeta function~$Z_{X,\chi}$ admits the factorization
\[
Z_{X,\chi}(s) = e^{p(s)} G_{X^\wedge,\chi}(s) G_\infty(s)^{-\dim(V) \chi_e^{\text{top}}(X)}\,\frac{\gammafunc\bigl(s-\tfrac12\bigr)^{n_p(\chi)}}{\gammafunc\bigl(s+\tfrac12\bigr)^{n_d(\chi)}}\, \mc P_{X,\chi}(s)
\]
on all of~$\CC$.
\end{thm}

In Theorem~\ref{thm:DP}, the function~$G_{X^\wedge,\chi}$ is an entire function with zeros only at the non-negative integers that encodes the contribution of the orbifold points, hence of the elliptic elements of~$\Gamma$. We refer to~\cite{Doll_Pohl} for precise formulas and emphasize here that $G_{X^\wedge,\chi}$ indeed encodes the contribution of each orbifold point separately. The numbers~$n_p(\chi)$ and~$n_d(\chi)$ are the sum of the singularity degrees at cusps and disk ends, respectively.

Therefore Theorem~\ref{thm:DP} provides not only a divisor formula for~$Z_{X,\chi}$ and hence a very precise relation between the zeros of~$Z_{X,\chi}$ and the $\chi$-twisted Laplace resonances of~$X$. It also shows that the contribution of each (coarse) geometric property of~$X$ (funnel ends, cusp ends, orbifold points) can be traced separately, therefore establishing that the geometry at infinity of~$X$ (in the presence of~$\chi$) has a great influence on the result. More precisely, the gamma factor~$\gammafunc\bigl(s-\tfrac12\bigr)^{n_p(\chi)}$ is completely determined by the $\chi$-twisted geometry at cusps, whereas the gamma factor~$\gammafunc\bigl(s+\tfrac12\bigr)^{n_d(\chi)}$ is completely determined by the $\chi$-twisted geometry at disk ends. The topological Euler characteristic~$\chi_e^{\text{top}}(X)$ only depends on the (untwisted) geometry at infinity (both cusps and disk ends), and hence the factor~$G_\infty(s)^{-\dim(V) \chi_e^{\text{top}}(X)}$ is determined by the untwisted geometry at infinity and, in a distinct way, the dimension of the twist. Contrary to this, the geometry at infinity does not influence the factor~$G_{X^\wedge,\chi}(s)$, which depends only on the orbifold singularities. Finally, the spectral information of~$(X,\chi)$, i.e., here the set of resonances, only influences the Hadamard product~$\mc P_{X,\chi}$. Of course, this product
depends on the geometry at infinity but in a different, and much more involved fashion than the other factors in the factorization formula in Theorem~\ref{thm:DP}.

\section{Selberg zeta functions with non-unitary twists}
\label{sec:nonunitary}

In view of the results surveyed in Section~\ref{sec:unitary} it is natural question to which extent they can be generalized to non-unitary twists. Throughout this section let $X=\Gamma\backslash\HH$ be any geometrically finite hyperbolic orbisurface and let $\chi\colon\Gamma\to\GL(V)$ be any finite-dimensional representation, also allowing non-unitary representations.

To ensure convergence of the infinite product~\eqref{def:szftwist} (used in the definition of the twisted Selberg zeta function) on some right half-space in~$\CC$, additional properties of the admissible representations~$\chi$ need to be requested that can be motivated by a geometric argument as follows.

Periodic geodesics can travel arbitrarily far into any cusp of~$X$ and hence wind around a cusp (end) arbitrarily often. As cusps correspond to (conjugacy classes of) parabolic elements in~$\Gamma$, the (conjugacy class of the) hyperbolic element corresponding to such a (primitive) periodic geodesic is of the form 
\[
 g_m \coloneqq g_2 p^m g_1
\]
with $p\in\Gamma$ being a parabolic element encoding the considered cusp, $m$ being the number of winds around the cusp, and $g_1,g_2\in\Gamma$ being suitable elements encoding the passage of the geodesic outside the cusp end. The contribution by~$\chi$ of this element in the infinite product~\eqref{def:szftwist} is
\[
 \chi(g_m) = \chi(g_2)\chi(p)^m\chi(g_1)\,.
\]
As $m$ may essentially be any integer (reflecting any number of windings into any direction around the cusp), an exponential contribution in~$m$ of the factor~$\chi(p)^m$ needs to be avoided to ensure convergence of~\eqref{def:szftwist}. For that reason we request for the representation~$\chi$ that for each parabolic element~$p\in\Gamma$, the eigenvalues of~$\chi(p)$ are of modulus~$1$. This property is known as \emph{non-expanding cusp monodromy}, for short NECM. As the following theorem states, the NECM representations are precisely those for which a reasonable twisted Selberg zeta function can be expected.

\begin{thm}[\cite{FP_szf}]\label{thm:FP}
If $\chi$ is NECM, then the infinite product~$Z_{X,\chi}$ in~\eqref{def:szftwist} converges for $\Rea s\gg 1$. If $\chi$ is not NECM, then this infinite product does not converge at all.
\end{thm}

Ongoing work with L.~Breitkopf indicates that a similar result holds for real hyperbolic spaces of higher dimension.

Regarding hyperbolic orbisurfaces, it is expected that Theorem~\ref{thm:DP} for the divisor can be generalized to Selberg zeta functions twisted by NECM representations. For generalizing the approach used in~\cite{GuZw95, DFP, DFP2, Doll_Pohl}, an important step is to understand the structure of Fourier-type expansions of functions satisfying the equivariance condition~\eqref{eq:equivariance}. As the statement of the result obtained in~\cite{FPR} is rather lengthy, we provide here an informal version only. We refer to~\cite{FPR} for a precise statement.

\begin{thm}[\cite{FPR}, informal version]\label{thm:FPR}
Let $A\in \GL(V)$ and let $f\colon\HH\to V$ be a vector-valued Laplace eigenfunction that satisfies the twist-periodicity condition 
\[
 f(z+1) = A f(z)
\]
for all $z\in\HH$. Then $f$ admits a Fourier expansion of the form 
\[
 f(z) = \sum_{n\in\ZZ} BJ^xB^{-1} \hat f_n(y,s) e^{2\pi i n x} \qquad (z=x+iy\in\HH)\,,
\]
where $J$ is the Jordan-like normal form of~$A$, and $B$ a transformation matrix to the Jordan-like normal form. The functions $\hat f_n$ for $n\not= 1$ are vector-valued. They essentially consist of modified Bessel functions of first and second kind and their derivatives. The zeroth Fourier coefficient is, as in the scalar case, formed in a different way, subject to resonance situations.
\end{thm}

In the case $A=1$, Theorem~\ref{thm:FPR} reduces to the classical Fourier expansion of periodic Laplace eigenfunctions. In the general case, the investigations can be done separately for each Jordan-like block. The size of the Jordan-like blocks determines the order of derivatives as well as the precise structure of the zeroth Fourier coefficient function. 

At least for establishing the meromorphic continuation of the infinite product and hence establishing a decent NECM-twisted Selberg zeta function, an alternative approach is available by transfer operator techniques. In a nutshell, a transfer operator is a weighted evolution operator that derives from a well-chosen discretization of the geodesic flow on the considered hyperbolic orbisurface. In a certain sense, a transfer operator serves as a mediator between the geodesic flow and the Laplace operator. A \emph{strict transfer operator approach} to the (twisted) Selberg zeta function refers to the existence of a discretization of the geodesic flow such that the $s$-weighted family of associated transfer operators~$(\mc L_{\chi,s})_{\Rea s \gg 1}$ represents the Selberg zeta function as a Fredholm determinant
\[
 Z_{X,\chi}(s) = \det\left( 1 - \mc L_{\chi, s}\right)\qquad (\Rea s \gg 1)
\]
and that $(\mc L_{\chi, s})_s$ admits a meromorphic continuation in~$s$ to all of~$\CC$. The latter property yields a meromorphic continuation of~$Z_{X,\chi}$ to all of~$\CC$.

\begin{thm}[\cite{FP_szf}]\label{thm:FP2}
If $X$ admits a strict transfer operator approach to the untwisted Selberg zeta function of~$X$, then it admits a strict transfer operator approach to the Selberg zeta function of~$X$ twisted by any NECM representation.
\end{thm}

As a consequence, Theorem~\ref{thm:FP2} guarantees the meromorphic continuation of the twisted Selberg zeta function if the untwisted Selberg zeta function can be handled by a strict transfer operator approach. The NECM property of representations is indeed crucial in the proof. 

Regarding the known realm of existence of such strict transfer operator approaches, the current state of art is~\cite{Pohl_Wabnitz, Wabnitz} in combination with~\cite{Pohl_symdyn}.

\begin{thm}[\cite{Pohl_Wabnitz, Wabnitz,Pohl_symdyn}, informal version]\label{thm:strictTO}
Strict transfer operator approaches exist for a huge family of geometrically finite noncompact hyperbolic orbisurfaces (with and without cusps).
\end{thm}

We refer to~\cite{Pohl_symdyn, Wabnitz} for a precise definition of the families of hyperbolic orbisurfaces for which Theorem~\ref{thm:strictTO} is valid. Ongoing work indicates that strict transfer operator approaches are expected to exist for \emph{all} geometrically finite noncompact hyperbolic orbisurfaces.

\section{Higher dimensional spaces}

\emph{A priori}, it is rather surprising that twisting representations sometimes allow to pass from spaces of low dimension to spaces of higher dimensions. Together with R.~Bruggeman and Y.~Choie, the author of this note developed such a situation in~\cite{ABC}, where the higher dimensional space is not even a locally symmetric space.

To provide a brief survey, let $G \coloneq \SL_2(\RR)$ and $\Gamma \coloneqq \SL_2(\ZZ)$. For certain specific reasons we do not descent to Fuchsian groups here. Further let $\Hei$ denote the $3$-dimensional continuous Heisenberg group and let $\Hei(\ZZ)$ be its subgroup defined by restricting all parameters to~$\ZZ$, which is a discrete group. The continuous Jacobi group of level~$1$ is the semidirect product $G^J \coloneqq \Hei \rtimes G$, and $\Gamma^J \coloneqq \Hei(\ZZ) \rtimes \Gamma$ is a discrete subgroup. We consider the product space $\HH\times\CC$ of the hyperbolic plane~$\HH$ and the complex plane~$\CC$, thus leaving the world of hyperbolic orbisurfaces. The Jacobi group~$G^J$ acts on $\HH\times\CC$ by fractional linear transformations in the $\HH$-component and by a certain skew product in the $\CC$-component. We refer to~\cite{ABC} for precise formulas.

\emph{Jacobi Maass cusp forms} for~$\Gamma^J$ are functions $\HH\times\CC\to\CC$ that are smooth eigenfunctions of certain differential operators on~$\HH\times\CC$, are of rapid decay towards the end of $\Gamma^J\backslash (\HH\times\CC)$ and are equivariant under a certain action of~$\Gamma^J$ that encodes an index and a weight. The precise definitions are involved, for which reason we omit them here and refer to~\cite{ABC} for details.

In \cite{Pit09}, it is shown that the spaces of Jacobi Maass forms of integral weight and positive integral index are isomorphic (as vector spaces) to certain spaces of vector-valued Maass forms of half-integral weight on~$\Gamma\backslash\HH$, thus certain spaces of vector-valued Laplace eigenfunctions on~$\Gamma\backslash\HH$. Extending the proof method of a theta decomposition used in~\cite{Pit09}, this isomorphism is extended in~\cite{ABC} to arbitrary real weight and any integer index.

Furthermore, by developing cohomological methodology in the vein of~\cite{BLZ15, BCD, tahA}, we establish in \cite{ABC} a notion of so-called \emph{period functions} for vector-valued Maass cusp forms for~$\Gamma$ of any real weight and any finite-dimensional unitary representation, and hence, by means of the isomorphism based on the theta decomposition, for Jacobi Maass cusp forms. Such period functions are closely related to eigenfunctions of the transfer operator families that made an appearance in Section~\ref{sec:nonunitary}. We refer to~\cite{ABC} for all details and close here with a final remark.

Even though the methodology and techniques used for the results surveyed in this section differ considerably from those of the previous sections, with a view to the usage of representations, a common underlying principle exists in the sense that representations, as long as they are sufficiently compatible or combinable with the other techniques and entities involved, are both very helpful in descending to easier-to-handle situations and to simultaneously define twisted objects of independent interest.

\printbibliography

\end{document}